\documentclass[11pt]{amsart}
\usepackage{amsmath, amssymb}
\usepackage{amsfonts}
\usepackage[arrow,matrix,curve,cmtip,ps]{xy}

\usepackage{amsthm}

\allowdisplaybreaks

\newtheorem{theorem}{Theorem}[section]

\newtheorem*{theorem*}{Theorem}
\theoremstyle{remark}

\newtheorem{definition}[theorem]{Definition}
\newtheorem{example}[theorem]{Example}

\numberwithin{equation}{section}

\newcommand{\Z}{\mathbb{Z}}
\newcommand{\Li}{\mathcal{L}}
\newcommand{\K}{\mathcal{K}}
\newcommand{\clspan}{\operatorname{\overline{\mathrm{span}}}}

\newcommand{\TX}{\mathcal{T}_X}
\newcommand{\OX}{\mathcal{O}_X}
\newcommand{\im}{\operatorname{im }}

\begin{document}
\title[Strong Shift Equivalence]{Strong Shift Equivalence in the $C^*$-algebraic Setting: Graphs and $C^*$-correspondences}

\author{Mark Tomforde 
}

\address{Department of Mathematics\\ University of Iowa\\
Iowa City\\ IA 52242-1419\\ USA}

\curraddr{Department of Mathematics\\ College of William \& Mary\\ P.O. Box 8795 \\
Williamsburg\\ VA 23187-8795\\ USA}

\email{tomforde@math.wm.edu}

\thanks{This research was supported by NSF Postdoctoral Fellowship DMS-0201960}

\date{\today}



\begin{abstract}

We discuss strong shift equivalence, which has been used to characterize conjugacy of edge shifts, and its application to $C^*$-algebras of graphs and Cuntz-Pimsner algebras.

\end{abstract}

\maketitle

\section{Background and Motivation: Edge Shifts of Finite Graphs}

A directed graph (hereafter simply called a graph) is a quadruple $E=(E^0,E^1,r,s)$ consisting of a countable set of vertices $E^0$, a countable set of edges $E^1$, and functions $r: E^1 \to E^0$ and $s : E^1 \to E^0$ identifying the range and source of each edge.  We say a graph is finite if $E^0$ and $E^1$ are finite sets.  If $E$ is a finite graph, then one may create a (two-sided) shift space $X_E$ with shift map $\sigma_E :X_E \to X_E$ defined by $$X_E := \{ (e_i)_{i \in \Z} : e_i \in E^1 \text{ for all } i \in \Z \text{ and }r(e_i) = s(e_{i+1}) \}$$ and $\sigma_E ((e_i)_{i \in \Z}) = (e_{i+1})_{i \in \Z}$.  The shift spaces that arise in this way are called \emph{edge shifts}.  We refer the reader to \cite{LM} for more about shift spaces and edge shifts.

For a given graph $E$, the \emph{vertex matrix} $A_E$ is the non-negative $E^0 \times E^0$ matrix defined by $$A_E (v,w) = \# \{ e \in E^1 : s(e)=v \text{ and } r(e)=w \}.$$

\begin{definition} Let $A$ and $B$ be finite, square, non-negative, integer matrices.  We say that $A$ and $B$ are \emph{elementary strong shift equivalent} if there exist non-negative, integer matrices $R$ and $S$ such that $A = RS$ and $B=SR$.  (Note that $A$ and $B$ need not have the same size.)
\end{definition}

Elementary strong shift equivalence is a relation that is reflexive and symmetric, but not transitive.  Therefore, we consider the equivalence relation it generates:  

\begin{definition}
We say that two finite, square, non-negative, integer matrices $A$ and $B$ are \emph{strong shift equivalent} if there exists a finite sequence $C_1, C_2, \ldots, C_n$ of finite, square, non-negative, integer matrices with $C_1 = A$, $C_n = B$, and $C_i$ elementary strong shift equivalent to $C_{i+1}$ for $i=1, \ldots n-1$.
\end{definition}

In 1973, R.~F.~Williams proved a remarkable theorem \cite{Wil}, which states that if $E$ and $F$ are finite graphs, then the edge shifts $(X_E, \sigma_E)$ and $(X_F, \sigma_F$ are conjugate (i.e.~there exists a homeomorphism $h : X_E \to X_F$ with $ \sigma_F \circ h = h \circ \sigma_E$) if and only if the edge matrices $A_E$ and $A_F$ are strong shift equivalent.  Thus the conjugacy class of an edge shift $(X_E, \sigma_E)$ is completely characterized by the strong shift equivalence class of its vertex matrix $A_E$.

\section{Strong Shift Equivalence and Graph $C^*$-algebras} \label{graph-sec}

If $E = (E^0,E^1,r,s)$ is a (not necessarily finite) graph, then the \emph{graph $C^*$-algebra} $C^*(E)$ is the universal $C^*$-algebra generated by a collection of mutually orthogonal projections $\{p_v : v \in E^0 \}$ together with a collection of partial isometries $\{ s_e : e \in E^1 \}$ with mutually orthogonal range projections that satisfy the \emph{Cuntz-Krieger relations}:
\begin{enumerate}
\item $s_e^*s_e=p_{r(e)}$ \text{ for all $e \in E^1$}
\item $s_es_e^* \leq p_{s(e)}$ \text{ for all $e \in E^1$}
\item $p_v = \sum_{\{e \in E^1 : s(e)=v \} } s_es_e^*$ \text{ for all $v \in E^0$
with $0 < | s^{-1}(v) | < \infty$}.
\end{enumerate}

We use the conventions established in \cite{KPRR, KPR,BPRS, FLR, RS, BHRS} for graph $C^*$-algebras.  We also refer the reader to \cite{Rae} for a more comprehensive treatment of graph $C^*$-algebra theory --- although we warn the reader that the direction of the arrows in \cite{Rae} is ``opposite" of what is used in \cite{KPRR, KPR, BPRS, FLR, RS, BHRS} and of what is used here.

\begin{definition}
A graph $E = (E^0,E^1,r,s)$ is said to be \emph{regular} if each vertex emits a finite and nonzero number of edges; i.e. $0 < |s^{-1}(v)| < \infty$ for all $v \in E^0$.  A (possibly infinite) matrix $A$ with 
 entries is said to be \emph{regular} if every row of $A$ contains a finite and nonzero number of positive entries.  Note that a graph $E$ is regular if and only if its vertex matrix $A_E$ is regular.  
\end{definition}

We now extend the definition of strong shift equivalent to regular matrices.

\begin{definition}  Let $A$ and $B$ be regular, square, non-negative, integer matrices.  We say $A$ and $B$ are \emph{elementary strong shift equivalent} if there exist non-negative, integer matrices $R$ and $S$ such that $A =RS$ and $B=SR$.  (Note that in order for $RS$ and $SR$ to be regular it is necessary that $R$ and $S$ be regular.)  We say that $A$ and $B$ are \emph{strong shift equivalent} if there exists a finite sequence $C_1, C_2, \ldots, C_n$ of regular, square, non-negative, integer matrices with $C_1 = A$, $C_n = B$, and $C_i$ elementary strong shift equivalent to $C_{i+1}$ for $i=1, \ldots n-1$.
\end{definition}

When $E$ is a finite graph with no sinks, the $C^*$-algebra $C^*(E)$ is intimately related to the edge shift $(X_E, \sigma_E)$.  We have already described how in this case the strong shift equivalence class of $A_E$ determines the conjugacy class of $(X_E, \sigma_E)$.  As we shall see, this has implications for the $C^*$-algebra $C^*(E)$.  In particular, if two graphs have vertex matrices that are strong shift equivalent, then their associated graph $C^*$-algebras are Morita equivalent.  Interestingly, this result holds also for infinite graphs that are regular.

\begin{theorem} \label{SSE-graphs}
Let $E$ and $F$ be regular graphs.  If the vertex matrices $A_E$ and $A_F$ are strong shift equivalent, then $C^*(E)$ and $C^*(F)$ are Morita equivalent.
\end{theorem}

This theorem was proven for Cuntz-Krieger algebras (which correspond to $C^*$-algebras of finite graphs with no sinks) by Cuntz and Krieger in \cite[Theorem~3.8]{CK}.  The theorem was proven for $C^*$-algebras of regular graphs by Bates in \cite[Theorem~5.2]{Bat} and by Drinen and Sieben in \cite[Proposition~7.2]{DS}.

$ $

\noindent \textbf{Sketch of Proof:}  Since Morita equivalence is an equivalence relation, it suffices to verify the claim when $A_E$ and $A_F$ are elementary strong shift equivalent.  Suppose that $A_E = RS$ and $A_F = SR$.  We may then form a bipartite graph $G_{R,S}$ that has vertices $E^0 \sqcup F^0$, and for each $v \in E^0$ and $w \in F^0$ there are $R(v,w)$ edges from $v$ to $w$ and $S(w,v)$ edges from $w$ to $v$.  Because $A_E = RS$, the paths of length two beginning in $E^0$ form a copy of $E$ in $G_{R,S}$, and because $A_F = SR$, the paths of length two beginning in $F^0$ form a copy of $F$ in $G_{R,S}$.  One can then use the Gauge-Invariant Uniqueness Theorem \cite[Theorem~2.1]{BPRS} to prove that $C^*(G_{R,S})$ contains subalgebras isomorphic to of $C^*(E)$ and $C^*(F)$, and furthermore one can show that these subalgebras are complementary full corners of $C^*(G_{R,S})$ determined by the projections $P=\sum_{v \in E^0} P_v$ and $Q=\sum_{v \in F^0} P_v$.  (We mention that if these sums are infinite, one can show that they converge to a projection in the multiplier algebra.)  Thus $C^*(E)$ and $C^*(F)$ are Morita equivalent.

\begin{example} \label{graph-ex}
The techniques of this proof are perhaps best illustrated with an example.  Suppose $E$ and $F$ are the graphs
$$\begin{matrix}
E & & \xymatrix{ v \ar[r]^b \ar@(dl,ul)^a & w \ar@(dr,ur)_c \\ } & & \qquad & &  F & &  \xymatrix{ x \ar[r]^e \ar@(dl,ul)^d & y \ar[r]^f &z \ar@(dr,ur)_g \\ }
\end{matrix}$$
Then we see that $A_E = \left( \begin{smallmatrix} 1 & 1 \\ 0 & 1 \end{smallmatrix} \right)$ and $A_F = \left( \begin{smallmatrix} 1 & 1 & 0 \\ 0 & 0 & 1 \\ 0 & 0 & 1 \end{smallmatrix} \right)$ are elementary strong shift equivalent by taking $R = \left( \begin{smallmatrix} 1 & 1 & 0 \\ 0 & 0 & 1 \end{smallmatrix} \right)$ and $S = \left( \begin{smallmatrix} 1 & 0 \\ 0 & 1 \\ 0 & 1 \end{smallmatrix} \right)$.  The bipartite graph $G_{R,S}$ is then equal to
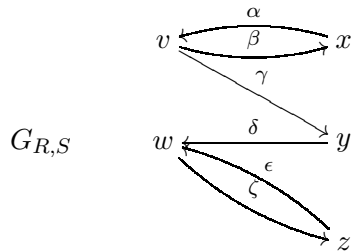
\begin{equation*}
\xymatrix{ & v \ar@/_/[rr]^\beta \ar[rrd]^\gamma & & x \ar@/_/[ll]_\alpha \\ G_{R,S} & w \ar@/_/[drr]^\zeta & & y \ar[ll]_\delta \\ & & & z \ar@/_/[llu]_\epsilon }
\end{equation*}
\end{example}
We see that the paths of length two in $G_{R,S}$ beginning in $E^0$ form a copy of $E$. (The edge $a$ corresponds to the path $\beta \alpha$, the edge $b$ corresponds to the path $\gamma \delta$, and the edge $c$ corresponds to the path $\zeta \epsilon$.)  Similarly, the paths of length two in $G_{R,S}$ beginning in $F^0$ form a copy of $F$.  Also, $C^*(E)$ and $C^*(F)$ are isomorphic to complementary full corners of $C^*(G_{R,S})$.  In fact, if $\{ s_e, p_v\}$ is a generating Cuntz-Krieger $E$-family for $C^*(E)$ and if $\{S_e, P_v \}$ is a generating Cuntz-Krieger $G$-family for $C^*(G_{R,S})$, then the $*$-homomorphism that identifies $C^*(E)$ with a full corner of $C^*(G_{R,S})$ maps
\begin{equation*}
p_v \mapsto P_v, \quad p_w \mapsto P_w, \quad s_a \mapsto S_\beta S_\alpha, \quad s_b \mapsto S_\gamma S_\delta, \quad \text{ and } \quad s_c \mapsto S_\zeta S_\epsilon.
\end{equation*}
Furthermore, this subalgebra is equal to the corner determined by $P := P_v + P_w$, and this corner is full since any hereditary subset of $G_{R,S}$ containing $E^0$ must contain all the vertices of $G_{R,S}$.  Similarly, one can see that $C^*(F)$ is isomorphic to the full corner determined by $Q = P_x + P_y  + P_z$, and these corners are complementary since $P + Q =1$.

We mention that there are counterexamples to Theorem~\ref{SSE-graphs} if one of the graphs contains a vertex that emits either no edges or infinitely many edges.  Thus the regularity of the graphs is necessary.

\section{Strong Shift Equivalence and Cuntz-Pimsner Algebras}

In this section we shall discuss a notion of strong shift equivalence for $C^*$-correspondences, and discuss how the strong shift equivalence class of a (essential, regular) $C^*$-correspondence determines the Morita equivalence class of the associated Cuntz-Pimsner algebra.  We will begin by giving some basic definitions and establishing notation.  Afterward we shall discuss how graph $C^*$-algebras can be realized as Cuntz-Pimsner algebras, and use this to motivate a definition of strong shift equivalence for $C^*$-correspondences.  We will then conclude with a description of how the proof of Theorem~\ref{SSE-graphs} can be generalized to the Cuntz-Pimsner setting.

\subsection{$C^*$-correspondences}

If $X$ is a right Hilbert $A$-module we let $\Li(X)$ denote the
$C^*$-algebra of adjointable operators on $X$, and we let $\K (X)$ denote
the closed two-sided ideal of compact operators given by $$\K (X) := \clspan
\{ \Theta_{\xi,\eta}^X : \xi, \eta \in X \}$$ where $\Theta_{\xi,\eta}^X$ is
defined by $\Theta_{\xi,\eta}^X (\zeta) := \xi \langle \eta, \zeta
\rangle_A$.  When no confusion arises we shall often omit the superscript and
write $\Theta_{\xi,\eta}$ in place of $\Theta_{\xi,\eta}^X$.

\begin{definition}
If $A$ and $B$ are $C^*$-algebras, then a \emph{$C^*$-correspondence from $A$ to $B$} is a right
Hilbert $B$-module $X$ together with a $*$-homomorphism $\phi_X : A \to
\Li(X)$.  We consider $\phi_X$ as giving a left action of $A$ on $X$ by setting
$a \cdot x := \phi_X(a) x$.  When $X$ is a $C^*$-correspondence from $A$ to $B$ we will sometimes write ${}_AX_B$ to keep track of the $C^*$-algebras.  If $A=B$ we refer to $X$ as a \emph{$C^*$-correspondence over $A$}.
\end{definition}

\begin{definition}
A $C^*$-correspondence $X$ from $A$ to $B$ is said to be \emph{essential} if $\clspan \{ \phi_X(a)x : a \in A \text{ and } x \in X \} = X$.  A $C^*$-correspondence is said to be \emph{regular} if $\phi_X$ is injective and $\phi_X(A) \subseteq \K (X)$.
\end{definition}

\begin{definition} \label{repn-def}
If $X$ is a $C^*$-correspondence over $A$, then a \emph{representation} of
$X$ into a $C^*$-algebra $B$ is a pair $(t, \pi)$ consisting of a linear map $t : X \to B$ and a
$*$-homomorphism $\pi : A \to B$ satisfying
\begin{enumerate}
\item[(i)] $t(\xi)^* t(\eta) = \pi(\langle \xi, \eta \rangle_X)$
\item[(ii)] $t(\phi_X(a)\xi) = \pi(a) t(\xi)$
\item[(iii)] $t(\xi a) = t(\xi) \pi(a)$
\end{enumerate}
for all $\xi, \eta \in X$ and $a \in A$.  We often write $(t, \pi) : (X,A) \to B$ in this situation.

If $(t, \pi) : (X,A) \to B$ is a representation of $X$ into
a $C^*$-algebra $B$, we let $C^*(t, \pi)$ denote the $C^*$-subalgebra of
$B$ generated by $t(X) \cup \pi(A)$.
\end{definition}

\begin{definition}[The Toeplitz Algebra of a $C^*$-correspondence]
Given a $C^*$-correspondence $X$ over a $C^*$-algebra $A$, there is a $C^*$-algebra $\TX$ and a representation $(\overline{t}_X, \overline{{\pi}}_X) :(X,A) \to \TX$ that is universal in the following sense: 
\begin{enumerate}
\item $\TX$ is generated as a $C^*$-algebra by $\overline{t}_X(X) \cup \overline{\pi}_X(A)$; and 
\item Given any representation $(t, \pi) :(X,A) \to B$ of $X$ into a $C^*$-algebra $B$, there exists a $*$-homomorphism of ${\rho}_{(t, \pi)} : \TX \to B$, such that $t={{\rho}_{(t, \pi)}}\circ{\overline{t}_X}$ and $\pi={{\rho}_{(t, \pi)}}\circ{\overline{\pi}}_X$.  
\end{enumerate}
The $C^*$-algebra $\TX$ and the representation $(\overline{t}_X, {\overline{\pi}}_X)$ exist (see \cite{FR}, for example) and are unique up to an obvious notion of isomorphism.  We call $\TX$ \emph{the Toeplitz algebra of the $C^*$-correspondence $X$}, and we call $(\overline{t}_X, {\overline{\pi}}_X)$ \emph{a universal representation} of $X$ in $\TX$.
\end{definition}

The Toeplitz algebra is a very natural $C^*$-algebra associated with a $C^*$-correspondence; however, in many practical situations it is too large.  The appropriate $C^*$-algebra to associate with a $C^*$-correspondence is called the Cuntz-Pimnser algebra.  It turns out that the Cuntz-Pimsner algebra is a quotient of the Toeplitz algebra, and it can be defined in this way.  However, we will instead define it in terms of its universal property, which involves coisometric representations.

\begin{definition}
For a representation $(t, \pi) : (X,A) \to B$ of $X$ into a $C^*$-algebra $B$ there
exists a $*$-homomorphism $\psi_t : \K (X) \to B$ with the property that
$$\psi_t (\Theta_{\xi,\eta}) = t(\xi) t(\eta)^*.$$  
See \cite[p.~202]{Pim}, \cite[Lemma~2.2]{KPW}, and \cite[Remark~1.7]{FR} for details on the existence
of this $*$-homomorphism.  (We warn the reader that our map $\psi_t$ is denoted by $\pi^{(1)}$ in much of the literature, and by $\rho^{(t,\pi)} = \rho^{(\psi,\pi)}$ in \cite{FR}.  We have chosen to use $\psi_t$ because the map depends only on $t$ and not on $\pi$.) 
\end{definition}

\begin{definition}
For an ideal $I$ in a $C^*$-algebra $A$ we define $$I^\perp := \{ a \in A :
ab=0 \text{ for all } b \in I \}.$$  If $X$ is a $C^*$-correspondence over
$A$, we define an ideal $J(X)$ of $A$ by $J(X) := \phi_X^{-1}(\K(X))$.  We also
define an ideal $J_X$ of $A$ by $$J_X := J(X) \cap (\ker \phi_X)^\perp.$$  Note
that $J_X = J(X)$ when $\phi_X$ is injective, and that $J_X$ is the maximal
ideal on which the restriction of $\phi$ is an injection into $\K(X)$.
\end{definition}

\begin{definition}
If $X$ is a $C^*$-correspondence over $A$, then a representation $(t, \pi) : (X,A) \to B$ of $X$ into a $C^*$-algebra $B$ is said to be \emph{coisometric} if $$\psi_t (\phi_X(a)) = \pi(a) \qquad \text{for all $a \in J_X$.}$$
\end{definition}

\begin{definition}[The Cuntz-Pimsner Algebra of a $C^*$-correspondence]
Given a $C^*$-correspondence $X$ over a $C^*$-algebra $A$, there is a $C^*$-algebra $\OX$ and a coisometric representation $(t_X, {\pi}_X) :(X,A) \to \OX$ that is universal in the following sense: 
\begin{enumerate}
\item $\OX$ is generated as a $C^*$-algebra by $t_X(X) \cup {\pi}_X(A)$; and 
\item Given any coisometric representation $(t, \pi) :(X,A) \to B$ of $X$ into a $C^*$-algebra $B$, there exists a $*$-homomorphism of ${\rho}_{(t, \pi)} : \OX \to B$, such that $t={{\rho}_{(t, \pi)}}\circ{t_X}$ and $\pi={{\rho}_{(t, \pi)}}\circ{\pi}_X$.  
\end{enumerate}
The $C^*$-algebra $\OX$ and the representation $(t_X, \pi_X)$ exist (see \cite[\S4]{Kat4}) and are unique up to an obvious notion of isomorphism.  We call $\OX$ \emph{the Cuntz-Pimsner algebra of the $C^*$-correspondence $X$}, and we call $(t_X, {\pi}_X)$ \emph{a universal coisometric representation} of $X$ in $\OX$.  The universal property of $\TX$ allows one to see that there is an epimorphism $\rho : \TX \to \OX$ and thus $\OX$ is isomorphic to a quotient of $\TX$.
\end{definition}

\subsection{Viewing Graph $C^*$-algebras as Cuntz-Pimsner algebras}

Given a graph $E = (E^0,E^1,r,s)$ one may define a $C^*$-correspondence $X(E)$ over $A:= C_0(E^0)$ by letting 
\begin{equation*}
X(E) := \{ x : E^1 \to \mathbb{C} : \text{ the function } v \mapsto \sum_{
\{f \in E^1: r(f) = v \} } |x(f)|^2 \text{ is in $C_0(E^0)$} \ \}.
\end{equation*}
and giving $X(E)$  the operations 
\begin{align*}
(x \cdot a)(f) &:= x(f) a(r(f)) \text{ \quad for $f \in E^1$} \\
\langle x, y \rangle_{X(E)}(v) &:= \sum_{ \{ f \in E^1: r(f) = v \} }\overline{x(f)}y(f) \text{ \quad for $v \in E^0$} \\
(a \cdot x)(f) &:= a(s(f)) x(f) \text{ \quad for $f \in E^1$}.
\end{align*}
We call $X(E)$ the \emph{graph $C^*$-correspondence} associated to $E$, and it is a fact that $\mathcal{O}_{X(E)} \cong C^*(E)$ \cite[Proposition~4.4]{FR}.  

In particular, if we let $P_v := \pi_{X(E)}(\delta_v)$ and $S_e := t_{X(E)}(\delta_e)$, where $\delta_v$ and $\delta_e$ denote point masses, then $\{P_v, S_e : v \in E^0, e \in E^1 \}$ is a collection of projections and partial isometries that satisfy the Cuntz-Krieger relations and generate $\mathcal{O}_{X(E)}$.  Furthermore, the graph $E$ is regular if and only if the $C^*$-correspondence $X(E)$ is regular.  Also, $X(E)$ is always essential.

Thus the graph $C^*$-algebra may be thought of as the Cuntz-Pimsner algebra associated to the graph $C^*$-correspondence.  We refer the reader to \cite[\S 3]{MT} for a more detailed discussion and analysis of graph $C^*$-correspondences.

We shall now use the graph $C^*$-correspondence to generalize the notion of strong shift equivalence to essential, regular $C^*$-correspondences.  Suppose that $E$ and $F$ are regular graphs and that there are non-negative, integer matrices $R$ and $S$ with $A_E = RS$ and $A_F = SR$.  Then $R$ is an $E^0 \times F^0$ matrix, and we may create a bipartite graph $G_R$ by defining $G_R^0 := E^0 \sqcup F^0$ and for $v \in E^0$ and $w \in F^0$ we draw $R(v,w)$ edges from $v$ to $w$.  For this graph we may construct a $C^*$-correspondence $X_R$ from $A:= C_0(E^0)$ to $B := C_0(F^0)$ by setting
\begin{equation*}
X_R := \{ x : G_R^1 \to \mathbb{C} : \text{ the function } v \mapsto \sum_{
\{f \in G_R^1: r(f) = v \} } |x(f)|^2 \text{ is in $C_0(F^0)$} \ \}.
\end{equation*}
and giving $X_R$  the operations 
\begin{align*}
(x \cdot b)(f) &:= x(f) b(r(f)) \text{ \quad for $f \in G^1_R$} \\
\langle x, y \rangle_{X_R}(w) &:= \sum_{ \{ f \in G_R^1: r(f) = w \} }\overline{x(f)}y(f) \text{ \quad for $w \in F^0$} \\
(a \cdot x)(f) &:= a(s(f)) x(f) \text{ \quad for $f \in G_R^1$}.
\end{align*}

In a similar way, we define $G_S$ and a $C^*$-correspondence $X_S$ from $B := C_0(F^0)$ to $A:= C_0(E^0)$.

\begin{example}
If $R$ and $S$ are the matrices in Example~\ref{graph-ex}, then $G_R$ and $G_S$ are the following graphs:
$$
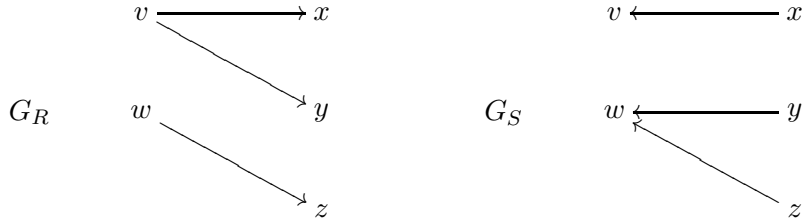
\begin{matrix}
 \xymatrix{ & v \ar[rr] \ar[rrd] & & x \\ G_R & w \ar[drr] & & y \\ & & & z }
 & \qquad & & \xymatrix{ & v & & x \ar[ll] \\ G_S & w & & y \ar[ll] \\ & & & z \ar[llu] }
\end{matrix}$$
\end{example}

The fact that $A_E = RS$ shows that any edge $e \in E^1$ corresponds to a path of length two $\alpha \beta$ with $\alpha \in G_R^1$ and $\beta \in G_S^1$.  Furthermore, this allows one to prove that $X(E) \cong X(G_R) \otimes_B X(G_S)$ via an isomorphism that sends $\delta_e \mapsto \delta_\alpha \otimes \delta_\beta$.  In a similar way we see that $A_F = SR$ implies that $X(F) \cong X(G_S) \otimes_A X(G_R)$.

This motivates the following definition.

\begin{definition}
Let $X$ be an essential, regular $C^*$-correspondence over $A$, and let $Y$ be an essential, regular $C^*$-correspondence over $B$.  We say $X$ and $Y$ are \emph{elementary strong shift equivalent} if there exists a $C^*$-correspondence $R$ from $A$ to $B$, and a $C^*$-correspondence $S$ from $B$ to $A$ such that $X =R \otimes_A S$ and $Y=S \otimes_B R$.  (We mention that if $X$ and $Y$ are essential and regular, then it follows from \cite[Corollary~3.11]{MPT} and \cite[Lemma~3.1.3]{MPT} that $R$ and $S$ may be chosen essential and regular.)  We say that $X$ and $Y$ are \emph{strong shift equivalent} if there exists a finite sequence $C_1$, $C_2$, \ldots $C_n$ of essential, regular $C^*$-correspondences with $C_1 = X$, $C_n = Y$, and $C_i$ elementary strong shift equivalent to $C_{i+1}$ for $i=1, \ldots n-1$.
\end{definition}

We then have the following generalization of Theorem~\ref{SSE-graphs}.

\begin{theorem} \label{SSE-correspondences}
Let $X$ be an essential, regular $C^*$-correspondence over $A$, and let $Y$ be an essential, regular $C^*$-correspondence over $B$.  If $X$ and $Y$ are strong shift equivalent, then $\OX$ and $\mathcal{O}_Y$ are Morita equivalent.
\end{theorem}

This result appears as Theorem~3.14 of \cite{MPT}, and a full proof can be found there.  Here we shall give a sketch of the proof with references to the appropriate lemmas of \cite{MPT}.

$ $

\noindent \textbf{Sketch of Proof:}  Since Morita equivalence is an equivalence relation, it suffices to verify the claim when $X$ and $Y$ are elementary strong shift equivalent.  Suppose that $X = R \otimes_A S$ and $Y = S \otimes_B R$.  In analogy with the proof of Theorem~\ref{SSE-graphs} we create a $C^*$-correspondence $Z$ over $A \oplus B$, called the \emph{bipartite inflation of  $S$ by $R$}.  We let $Z := S \oplus R$, and give $Z$ the structure of a right Hilbert $A \oplus B$-module via 
\begin{align*}
(s,r) \cdot (a,b) := (sa, rb) \quad \text{ and } \quad \langle (r_1, s_1), (r_2, s_2) \rangle := (\langle r_1, s_1 \rangle, \langle r_2, s_2 \rangle )
\end{align*}
and we make $Z$ into a $C^*$-correspondence over $A \oplus B$ by defining the left action as $$\phi_{A \oplus B} (a,b)(s,r) := (\phi_B(b)s, \phi_A(a)r).$$

The way to verify the claim is then to show that $\OX$ and $\mathcal{O}_Y$ are isomorphic to complementary full corners of $\mathcal{O}_Z$.  We shall outline how this is done:

Let $(t_Z, \pi_Z) : (Z, A \oplus B) \to \mathcal{O}_Z$ be a universal coisometric representation of $Z$ into $\mathcal{O}_Z$.  We then define a representation $(t,\pi) : (X,A) \to \mathcal{O}_Z$ by setting $$t (r \otimes s) := t_Z(0,r) t_Z(s,0) \quad \text{ and } \pi(a) := \pi_Z(a,0).$$  It is proven in \cite[Lemma~3.8]{MPT} that $(t, \pi)$ is a coisometric representation, and therefore induces a homomorphism $\rho_{(t,\pi)} : \mathcal{O}_X \to \mathcal{O}_Z$.  An application of the Gauge-Invariant Uniqueness Theorem allows one to show that $\rho_{(t,\pi)}$ is injective and $\mathcal{O}_X$ is isomorphic to the subalgebra $C^*(t, \pi) := \im \rho_{(t,\pi)}$ of $\mathcal{O}_Z$.  Similarly, one can define $(t', \pi') : (Y,B) \to \mathcal{O}_Z$ and it is shown in \cite[Lemma~3.9]{MPT} that $(t', \pi')$ is a coisometric representation and that the induced homomorphism $\rho_{(t', \pi')} : \mathcal{O}_Y \to \mathcal{O}_Z$ is injective, so that $\mathcal{O}_Y$ is isomorphic to the subalgebra $C^*(t', \pi') := \im \rho_{(t',\pi')}$ of $\mathcal{O}_Z$.

Finally, we let $\{e_\lambda \}_{\lambda \in \Lambda}$ be an approximate unit for $A$, and we let $\{f_\lambda \}_{\lambda \in \Lambda}$ be an approximate unit for $B$.  Since $X$ and $Y$ are essential, it follows that $Z$ is essential, and thus $P := \lim_\lambda \pi_Z(e_\lambda, 0)$ and $Q := \lim_\lambda \pi(0, f_\lambda)$ converge to projections in the multiplier algebra of $\mathcal{O}_Z$.  Furthermore, the following hold:
\begin{itemize}
\item The fact that $X$ and $Y$ are essential and regular implies that $P \mathcal{O}_Z P = C^*(t, \pi)$ and $Q \mathcal{O}_Z Q = C^*(t', \pi')$.  (See the proof of \cite[Theorem~3.14]{MPT}.)
\item The fact that $X$ and $Y$ are regular implies that the corners $P \mathcal{O}_Z P$ and $Q \mathcal{O}_Z Q$ are full.  (See the proof of \cite[Theorem~3.14]{MPT}.)
\item $P+Q=1$ (see \cite[Lemma~3.12]{MPT}).
\end{itemize}
Thus $\mathcal{O}_X$ and $\mathcal{O}_Y$ are isomorphic to complementary full corners of $\mathcal{O}_Z$.

\section{Concluding Remarks}

\subsection{Non-essential $C^*$-correspondences}

In Theorem~\ref{SSE-correspondences} we required that the $C^*$-correspondences in question be essential and regular.  (Also, since all graph $C^*$-correspondences are essential, the $C^*$-correspondences covered by Theorem~\ref{SSE-graphs} are both essential and regular.)  As pointed out at the end of \S \ref{graph-sec}, the regularity condition is a necessary hypothesis, and there are known counterexamples if it is removed.  On the other hand, it is currently unknown whether the essential condition is necessary in Theorem~\ref{SSE-correspondences}.  Nonetheless, it has been shown in \cite{MPT} that the condition that the $C^*$-correspondence is essential may be replaced by the condition that the underlying $C^*$-algebra is unital; that is:

\begin{theorem}[Theorem~4.3 of \cite{MPT}]
Let $X$ be a regular $C^*$-correspondence over a $C^*$-algebra $A$, and let $Y$ be a regular $C^*$-correspondence over a $C^*$-algebra $B$.  Suppose that either $X$ is essential or $A$ is unital.  Also suppose that either $Y$ is essential or $B$ is unital.  If $X$ is elementary strong shift equivalent to $Y$, then $\mathcal{O}_X$ is Morita equivalent to $\mathcal{O}_Y$.
\end{theorem}

It is currently unknown whether the condition that the underlying $C^*$-algebra is unital can be removed (or weakened) in the above theorem.

\subsection{Morita equivalence at other levels}

Suppose that $X$ is an essential, regular $C^*$-correspondence over $A$ and that $Y$ is an essential, regular $C^*$-correspondence over $B$.  There are three natural questions that one can ask:

$ $

\noindent \textsc{Question 1:}  If $X$ and $Y$ are elementary strong shift equivalent, then is it necessarily the case that $X$ and $Y$ are Morita equivalent as $C^*$-correspondences (as defined in \cite{MS2})?

\smallskip

\noindent \textsc{Question 2:}  If $X$ and $Y$ are elementary strong shift equivalent, then is it necessarily the case that the Toeplitz algebras $\mathcal{T}_X$ and $\mathcal{T}_Y$ are Morita equivalent?

\smallskip

\noindent \textsc{Question 3:}  If $X$ and $Y$ are elementary strong shift equivalent, then is it necessarily the case that the Cuntz-Pimsner algebras $\mathcal{O}_X$ and $\mathcal{O}_Y$ are Morita equivalent?

\smallskip

\smallskip

We have seen that Theorem~\ref{SSE-correspondences} provides an affirmative answer to Question~3.  In addition, we mention that an affirmative answer to Question~1 implies an affirmative answer to Question~2.  Furthermore, since the Cuntz-Pimsner algebra is a quotient of the Toeplitz algebra, we see that if $\mathcal{T}_X$ is Morita equivalent to $\mathcal{T}_Y$, and if the Morita equivalence takes the appropriate ideal in $\mathcal{T}_X$ to the appropriate ideal in $\mathcal{T}_Y$ via the Rieffel correspondence, then $\mathcal{O}_X$ is Morita equivalent to $\mathcal{O}_Y$.

Originally it was hoped that one could prove an affirmative answer to Question~1, and then use this result to obtain affirmative answers to Question~2 and Question~3 as described in the previous paragraph.  Surprisingly, however, it has been shown that this will not work.  In the example described in \cite[\S5.2]{MPT} it is  shown that there are $C^*$-correspondences $X$ and $Y$ (in fact, both can be chosen to be graph $C^*$-correspondences) with the property that the Toeplitz algebra $\mathcal{T}_X$ is not Morita equivalent to the Toeplitz algebra $\mathcal{T}_Y$.  It, of course, also follows that the $C^*$-correspondences $X$ and $Y$ are not Morita equivalent.  Consequently this example provides a counterexample to Question~1 and Question~2, and it shows that strong shift equivalence of $C^*$-correspondences implies Morita equivalence at the level of Cuntz-Pimsner algebras, but not at the level of $C^*$-correspondences or at the level of Toeplitz algebras.

\end{document}